\documentclass[aos,MSNbibl,seceqn,citesort,dvips]{arximspdf}
\usepackage{epic}
\usepackage{graphicx}
\doi{10.1214/11-AOS915}
\volume{39}
\issue{5}
\pubyear{2011}
\firstpage{2716}
\lastpage{2739}

\makeatletter

\newtheorem{prop}{Proposition}[section]

\newproclaim{defn}{Definition}

\newcommand{\para}{{\Vert}}
\newcommand{\splitsterm}{S_{\mathrm{term}}}
\newcommand{\spnpos}{\mathrm{span}_{+}}
\newcommand{\R}{\mathbb{R}}
\newcommand{\treespace}{\mathcal{T}_O}
\newcommand{\treespaceint}{\mathcal{T}_{O,\mathrm{int}}}

\makeatother

\begin{document}
\begin{frontmatter}

\title{Principal components analysis in the space of phylogenetic trees}
\runtitle{Principal components for phylogenies}

\begin{aug}
\author[A]{\fnms{Tom M. W.} \snm{Nye}\corref{}\ead[label=e1]{tom.nye@ncl.ac.uk}\ead[label=u1,url]{http://www.mas.ncl.ac.uk/\textasciitilde ntmwn}}
\runauthor{T. M. W. Nye}
\affiliation{Newcastle University}
\address[A]{School of Mathematics and Statistics\\
Newcastle University\\
Newcastle upon Tyne \\
Tyne and Wear \\
NE1 7RU \\
United Kingdom\\
\printead{e1}\\
\printead{u1}} 
\end{aug}

\received{\smonth{7} \syear{2010}}
\revised{\smonth{7} \syear{2011}}

%
\begin{abstract}
Phylogenetic analysis of DNA or other data commonly gives rise to a
collection or sample of inferred evolutionary trees. Principal
Components Analysis (PCA) cannot be applied directly to collections of
trees since the space of evolutionary trees on a fixed set of taxa is
not a vector space. This paper describes a novel geometrical approach
to PCA in tree-space that constructs the first principal path in an
analogous way to standard linear Euclidean PCA. Given a data set of
phylogenetic trees, a geodesic principal path is sought that maximizes
the variance of the data under a form of projection onto the path. Due
to the high dimensionality of tree-space and the nonlinear nature of
this problem, the computational complexity is potentially very high, so
approximate optimization algorithms are used to search for the optimal
path. Principal paths identified in this way reveal and quantify the
main sources of variation in the original collection of trees in terms
of both topology and branch lengths. The approach is illustrated by
application to simulated sets of trees and to a set of gene trees from
metazoan (animal) species.
\end{abstract}

%

\end{frontmatter}

\section*{Introduction}\label{intro}
Inference of evolutionary or \textit{phylogenetic} trees is a
fundamental task in many areas of biology, and tree estimation has
developed over several decades into a mature statistical field
\cite{fel04}. On a phylogenetic tree, leaves correspond to existing
observed taxa, internal vertices correspond to ancestral taxa, and
branch lengths represent the degree of evolutionary divergence between
taxa. A~phylogenetic tree representing the division and divergence of
different species is called a \textit{species tree}. However,
individual regions of DNA can evolve according to trees that differ
from the underlying species tree, and an inferred phylogenetic tree
from a particular gene or DNA region is called a \textit{gene tree}.
Gene trees can differ from the species tree for several reasons: random
variation in the process of DNA letter substitution; population effects
by which the evolutionary course of an individual gene does not match
that of the species as a whole \cite{don95}; and even relatively rare
events whereby genetic material is exchanged between species in a
nontree-like manner \cite{doo99}. Phylogenetic analysis of a number of
different genes in a fixed set of species therefore generally gives
rise to a collection of alternative phylogenetic trees.
Collections\vadjust{\goodbreak} of alternative phylogenetic trees also
arise from inferential methods that involve simulation: bootstrap
replication and MCMC sampling from Bayesian posteriors are widely used
in the construction of phylogenetic estimates. Given such a collection
of alternative trees, whether gene trees or a simulated sample,
identifying differences and quantifying variation is a difficult
problem, since we might potentially have several hundred trees on
thousands of species. Standard multivariate statistical methods such as
clustering \cite{sto02,chak10,nye08} and Multi-Dimensional Scaling
(MDS) \cite{hill05,chak10} have been used to address this problem.
Principal Components Analysis (PCA), in contrast, cannot be applied
directly since the space of phylogenetic trees on a fixed set of
species is not a Euclidean vector space. This paper describes a
geometric approach to PCA for sets of alternative phylogenetic trees.
The aim is to identify which tree features are most variable within a
given set of trees and to quantify this variation---just as the first
few components in regular PCA pick out the most variable features of a
Euclidean data set. Although PCA has been used to analyze different
phylogenetic data previously (such as distance matrix data), the method
presented here is the first to work intrinsically within the space of
phylogenetic trees. The approach relies to a large extent on existing
mathematical tools, and the main novel contribution comes from
combining those elements into a computationally feasible method.

\begin{figure}

\includegraphics{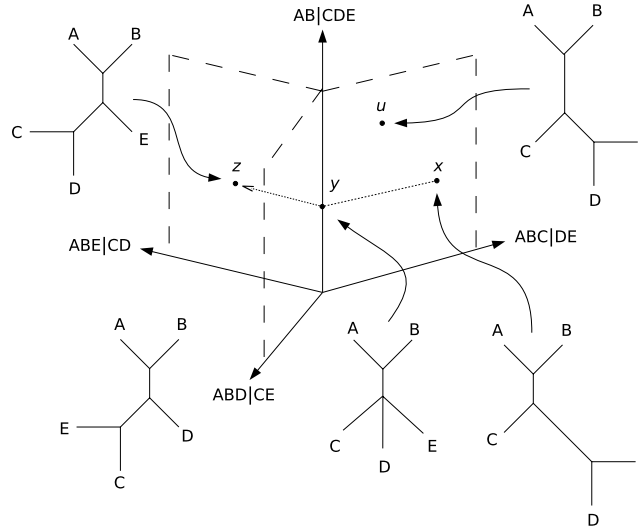}

\caption{A schematic view of a region of tree-space on
five taxa: points in space correspond to unrooted trees. Trees with the
same topology all lie in the same quadrant of tree-space (trees $x,u$,
e.g.). Different quadrants are joined along their edges. Tree $x$ can
be continuously deformed into tree $z$ by shrinking down the branch
$\mathit{DE}$ via tree $y$ and replacing it with the branch $\mathit{CD}$. It follows
that tree $z$ is obtained from $x$ via nearest neighbor interchange
(NNI) of the split $\mathit{ABC}|\mathit{DE}$ into split
$\mathit{ABE}|\mathit{CD}$. However, at $y$ another NNI move is
possible: $\mathit{ABC}|\mathit{DE}$ could be replaced by
$\mathit{ABD}|\mathit{CE}$, corresponding to the lower left quadrant.}
\label{figstitch}
\vspace*{-3pt}
\end{figure}

A key feature of our approach is the incorporation of both topological
and geometrical information from the trees under analysis, via the
so-called \textit{geodesic metric} on the space of trees \cite
{bill01,kup08,owen09}.
\textit{Topological} information refers to the exact pattern of branching
within a tree, while \textit{geometrical} information refers to the
distances between taxa induced by branch lengths on the tree.
Topological features are generally more straightforward to characterize
in a set of alternative trees, by counting the proportion of trees
containing a given feature.
For example, bootstrap replicate data sets are often represented by a
single ``consensus'' tree annotated with a percentage support for each
clade within the tree \cite{fel85}.
However, the geometry and topology of evolutionary trees are intimately
related: we can continuously change the topology of a tree by shrinking
down the length of any internal branch and expanding out an alternative
branch in its place, as shown in Figure \ref{figstitch}.
%
%
Recent authors \cite{kup08} have stressed the importance of using
geometrical information to draw comparisons between trees on account of
the interdependence of tree geometry and topology, and due to the
increased distinguishability obtained by using continuous rather than
discrete metrics.
Moreover, tree geometry plays an important role in inference: it has
been shown that long branches tend to ``attract'' each other, leading to
mistakes in the topology of inferred trees \cite{hue93}.

Taking a set of alternative phylogenetic trees on some fixed set of
taxa as input, our approach identifies a path $L$ through tree-space
that can be thought of analogously to the first principal component in
regular Euclidean PCA. The path consists of a smoothly changing tree
structure in which certain branches expand or shrink. Alternative
topologies emerge when internal branches are shrunk\vadjust{\goodbreak}
to have zero length and are then replaced with topologically distinct
branches. The path $L$ is constructed in such a way that the changing
features---both in terms of topology and geometry---correspond to the
most variable features within the data set. Just as for regular PCA,
$L$ also captures correlations in the data set: features that tend to
occur together in the data are also represented together on $L$. A
quantitative measure of variability can be assigned to $L$, in analogy
with the proportion of variance contributed by each component in
regular PCA. Unlike Euclidean vector spaces, there is no inner product
on tree-space, and so the analysis cannot be extended in a
straightforward manner to provide higher order principal paths by
working orthogonally to $L$. Further discussion is given in
Section~\ref{secconc}.\vadjust{\goodbreak}

Our approach---which we will refer to as $\Phi$PCA
(for ``phylogenetic''
PCA)---is motivated by geometrical analogy with regular vector space PCA.
Construction of the first principal component in a Euclidean vector
space can be thought of as follows:

{\renewcommand\thelonglist{(\arabic{longlist})}
\renewcommand\labellonglist{\thelonglist}
\begin{longlist}
\item Given a set of vectors $x_1, \ldots, x_n$ identify the centroid
$\bar{x}$.
\item For a fixed line $L$ through $\bar{x}$ take the orthogonal
projection of the points $x_1, \ldots, x_n$ onto $L$.
\item Identify the line that maximizes the variance of the projected
points along~$L$, or, equivalently, which minimizes the sum of squared
orthogonal distances of the points from the line.
\end{longlist}}

\noindent For a Euclidean vector space, these steps can be re-expressed and
solved in terms of simple linear algebra.
However, the space of phylogenetic trees on a~fixed set of taxa is not
a Euclidean vector space, so these steps cannot be applied directly in
the same way to sets of alternative phylogenies.
Tree-space can be equipped with various metrics that allow geometry to
be performed, and for reasons described below, we use the geodesic
metric \cite{bill01}.
$\Phi$PCA then follows a similar set of steps to those above, but working
with the geodesic metric, $d(\cdot, \cdot)$.
In step (2), the lines $L$ become paths in tree-space with the property
that, for any pair of points on the path, the path coincides with the
geodesic between the points.
Trees $x_1, \ldots, x_n$ are ``projected'' onto each path by finding
points $y_i$ on the path that minimize the distance $d(x_i,y_i)$ for
$i=1, \ldots, n$.
Pythagoras' theorem does not hold in tree-space, so in the analog of
step (3), paths which maximize the variance can be different from paths
which minimize the sum of squared distances.
We consider searching for both types of paths.
Step (3) is potentially excessively computationally demanding, and so
we describe (i) a greedy algorithm for constructing optimal
paths and (ii) a Monte Carlo optimization approach.
The methods we propose for steps (2) and (3) form the novel
contribution of this paper.
``Projection'' of points onto a geodesic path $L$ in step (2) is
relatively simple to perform using existing methods for computing the
geodesic metric, but a detailed algorithm has not been given previously.
Searching over the set of possible paths is more technically demanding.
Consideration of this particular problem and the solutions we present
appear to be entirely novel.

The development of $\Phi$PCA has been influenced by a recent paper by
Wang and Marron \cite{wang07}.
Wang and Marron addressed a similar problem, developing a form of PCA
for data sets with a tree-like structure.
In a~second paper \cite{ayd09}, they applied their method to sets of
trees obtained from medical imaging data.
In particular, their reformulation of PCA in terms of the geometrical
steps specified above motivated the corresponding steps in $\Phi$PCA.
Other authors have also developed analogs of PCA in nonstandard
geometries \cite{has89,flet04}, and Wang and Marron give an excellent
overview of this area of research \cite{wang07}.
However, it must be stressed that the method of Wang and Marron does
not apply to sets of phylogenetic trees, and that $\Phi$PCA is not simply
a reworking of their approach.
On account of the ostensible similarities between the approaches, we
devote a section to explaining the relationship between them later in
the paper.

The remainder of the paper is structured in the following way.
We first describe the geometry of tree-space and set up necessary
notation and mathematical background.
Section \ref{secmethods} contains a description of the $\Phi$PCA approach
and proofs of its properties.
We then explain more fully the relationship to the work of Wang and
Marron, before evaluating $\Phi$PCA on simulated sets of trees and a real
set of gene trees from metazoan species.


\section{Background: The geometry of tree-space}\label{secbackground}

\subsection{Splits and vector representation of trees}\label{secsplits}

We will work throughout with a fixed set of taxa $O=\{o_1, \ldots,
o_m\} $ and the set of unrooted phylogenetic trees $\treespace$ on $O$.
Given a tree $x\in\treespace$, cutting any branch on $x$ partitions the
taxa into two unordered nonoverlapping sets. Such a
partition is called a \textit{split}, and splits are usually denoted
$X|X^c$ where $X\subset O$ and $c$ denotes the complement in $O$. There
are $M=2^{m-1}-1$ possible (nonempty) splits of the set $O$, and the
set of these is denoted $S$. It is crucial to note that arbitrary sets
of splits do not generally correspond to valid tree topologies---a
compatibility condition must be satisfied. For example, if
$O=\{A,B,C,D,E\}$, then the two splits $\{A,B\}|\{ C,D,E\}$ and
$\{A,C\}|\{B,D,E\}$ cannot both be represented on the same tree.

Any tree $x\in\treespace$ can be regarded as a weighted set of
compatible splits, where the weight assigned to each split is given by
the length of the corresponding branch on $x$.
We only consider trees with positive branch lengths.
We write $T_x$ to denote the set of splits in $x$, and encapsulate the
branch lengths via a function $\lambda_x\dvtx S\rightarrow\R^+$
defined by
\[
\lambda_x(p) = \cases{
\mbox{branch length associated with $p$}, &\quad if $p\in T_x$,\cr
\mbox{zero}, &\quad otherwise.}
\]
Tree-space $\treespace$ can then be embedded in $\R^M$ in the
following way.
Take the standard basis of $\R^M$ and associate each split $p\in S$
with a different basis vector $\mathbf{e}_p$.
Any tree $x\in\treespace$ can then be associated uniquely with the vector
%
%
\begin{equation}\label{equembed}
\bolds\lambda_x = \sum_{p\in S} \lambda_x(p)\mathbf{e}_p.
\end{equation}
In fact, it is convenient to abuse notation slightly and write
$\mathbf{p}$ for the basis vector~$\mathbf{e}_p$, identifying each
split directly with the corresponding vector in~$\R^M$.
Equation (\ref{equembed}) essentially associates every tree $x$ with
a vector of branch lengths, but due to the compatibility relations
between splits, not every such vector corresponds to a tree.
In fact, each tree contains at most $2m-3$ splits, so as the number of
taxa $m$ increases, $2m-3\ll M$ and $\treespace$ becomes an
increasingly sparse subset of $\R^M$.

Since a collection $x_1, \ldots, x_n$ of trees can be regarded as a set
of vectors $\bolds\lambda_{x_1},\ldots, \bolds\lambda_{x_n}$,
why not just perform PCA on these vectors?
In general, the principal components obtained in this way will not
correspond to valid trees, and interpretation of the principal
components becomes impossible.
A form of PCA which operates intrinsically within $\treespace$ and
which produces interpretable ``components'' is required.

\subsection{Decomposition of tree-space by topology}

The geometry of $\treespace$ was first comprehensively studied in a
paper by Billera et al. \cite{bill01}, which included the definition
and proof of existence of geodesics.
Their description of $\treespace$ amounts to a decomposition into a set
of overlapping component pieces, each piece corresponding to a
different tree topology.
In this section we recall aspects of this decomposition which are
central to $\Phi$PCA, most importantly for the definition of geodesics
on~$\treespace$.

The decomposition is easiest to understand by identifying $\treespace$
with its image under the embedding in $\R^M$.
Every tree in $\treespace$ contains the set of splits corresponding to
terminal edges (those that end in a leaf), denoted $\splitsterm\subset
S$.
Since every tree contains every terminal split
\[
\treespace\cong\spnpos\{ \mathbf{p} \dvtx p\in\splitsterm\}\times
\treespaceint,
\]
where $\spnpos$ denotes the span of vectors with nonnegative weights,
and~$\treespaceint$ is the part of tree-space corresponding to
internal splits.
Next consider a~single unrooted tree $x$ which is fully resolved, by
which we mean every internal vertex has exactly 3 neighbors.
Let $t$ denote the topology of $x$ or, more precisely, the set of
nonterminal splits $t = T_x\setminus\splitsterm$.
Since $x$ is fully resolved, it has $m-3$ internal edges, so $t$
contains $m-3$ splits.
The internal branch lengths of any tree with topology $t$ are
determined by a point in $Q_t=\spnpos\{ \mathbf{p} \dvtx p\in t \}$.
We call $Q_t$ the \textit{topological orthant} containing $x$, and it is
isomorphic to the positive orthant of $\R^{m-3}$.
The faces of the orthant~$Q_t$ correspond to trees that have some zero
length branches.
Such trees are not fully resolved or, in other words, some internal
vertices have more than~3 neighbors.
This structure is illustrated in Figure \ref{figstitch}.

Tree-space $\treespace$ is formed from the union of the orthants $Q_t$
over all possible fully resolved topologies $t$:
\[
\treespace= \spnpos\{ \mathbf{p} \dvtx p\in\splitsterm\} \times
\mathop{\bigcup_{\mathrm{resolved}}}_{\mathrm{topologies}\ t}Q_t.
\]
The individual orthants $Q_t$ are stitched together along their faces,
since each unresolved tree occurs on the face of more than one
orthant.\vadjust{\goodbreak}
To understand how the orthants are stitched together in more detail,
consider a~point on the face of an orthant $Q_t$ at which a single
branch length corresponding to a split~$p$ has been collapsed to zero.
As illustrated in Figure \ref{figstitch}, there are two ways in which
this branch can be replaced with an alternative, thereby obtaining a
fully resolved tree with a different topology.
Each \mbox{$(m-4)$}-dimensional face of $Q_t$ is therefore identified with
corresponding faces in two other orthants $Q_{t'}$ and $Q_{t''}$.
The operation illustrated in Figure~\ref{figstitch} is called \textit
{Nearest Neighbor Interchange} (or NNI); we say that topologies~$t'$
and $t''$ are obtained by NNI of the split $p$ within $t$.
Faces of $Q_t$ with co-dimension greater than $1$ will be contained in
more than two other orthants.
Later, we will need to deal with paths in $\treespace$ between such
faces and so we need to extend the definition of NNI (which is usually
taken as a relationship between strictly binary trees).
Given a split $p$ in a fixed tree, there are two or more subtrees
hanging from each end of the associated edge in the tree.
An \textit{extended} NNI move (or XNNI) consists of swapping a subtree
from one end of the branch with a~subtree from the opposite end.
This operation removes split $p$ from the tree and replaces it with an
incompatible split $p'$.
On a~binary tree this definition coincides with the standard definition
of NNI \cite{all01}, and XNNI includes all NNI moves.

\subsection{Geodesics and the geodesic metric}\label{secgeodesicmetric}

$\treespace$ can be equipped with metrics via the embedding into $\R^M$
described above.
In particular, $L_2$ norm on~$\R^M$ defines a metric: $d_2(x,y) =
|\bolds\lambda_x - \bolds\lambda_y|_2$.
However, such metrics are not intrinsic to tree-space.
For example, when $x$ and $y$ have different topologies, $d_2$
corresponds to the length of a straight line segment joining $x$ to $y$
through $\R^M$, but this line contains points outside the image of
$\treespace$ under the embedding.

Billera et al. \cite{bill01} proved the existence of a metric that
locally resembles the $L_2$ metric, but which is intrinsic to
$\treespace$ independent of the embedding in $\R^M$.
This metric is called the \textit{geodesic} metric $d$, and it is the
canonical metric for $\Phi$PCA due to its intrinsic nature.
It is defined as follows.
For two trees $x$ and $y$ with the same topology,
$d(x,y)=d_2(x,y)$.
When $x$ and~$y$ have different topologies, $d(x,y)$ is defined as the
length of the shortest continuous path joining $x$ to $y$ in
$\treespace
$ which consists of a series of straight line segments through any
feasible sequence of topological orthants.
The length of such a path is defined to be the sum of the Euclidean
lengths of each constituent line segment.
The shortest such path joining $x$ to $y$ is called the \textit{geodesic}
between $x$ and $y$.
The proof that geodesics exist between points in~$\treespace$ and that
geodesics define a valid metric is given in \cite{bill01}.
As part of the proof, Billera et al. \cite{bill01} showed that
tree-space is CAT$(0)$ \cite{gro87}.
This means that triangles in $\treespace$ are ``skinny'' in comparison to
triangles in the Euclidean plane.
More formally, given points $x,y,z\in\treespace$, consider the triangle
between points $x',y',z'$ in the Euclidean plane with the same edge
lengths, so that $d(x,y)=d(x',y')$, etc.
If $\gamma(t)$ is the path-length parameterized geodesic between $x$
and $y$ and $\gamma'(t)$ the corresponding geodesic in the Euclidean
plane, then
$d(z,\gamma(t))\leq d(z',\gamma'(t))$ for all points $\gamma(t)$
between $x$ and $y$.

%
%
\begin{figure}
\setlength{\unitlength}{0.36cm}
\begin{center}
\begin{picture}(20,18)

\put(10.0,9.0){\vector(0,1){8}}
\put(10.0,9.0){\vector(0,-1){8}}
\put(10.0,9.0){\vector(1,0){8}}
\put(10.0,9.0){\vector(-1,0){8}}
\put(8.2,17.3){$DE|ABC$}
\put(17.3,8.0){$AB|CDE$}
\put(8.2,0.0){$BE|ACD$}
\put(0.5,8.0){$AC|BDE$}

\matrixput(10.0,9.0)(0.25,0){26}(0,-0.25){26}{\line(1,-1){1.0}}%

\put(2.0,16.0){\line(-1,1){1}}
\put(2.0,16.0){\line(1,1){1}}
\put(2.0,15.0){\line(1,0){1}}
\put(2.0,15.0){\line(0,1){1}}
\put(3.0,15.0){\line(1,1){1}}
\put(3.0,15.0){\line(1,-1){1}}
\put(2.0,15.0){\line(-1,-1){1}}
\put(0.5,17.1){$A$}
\put(3.0,17.1){$C$}
\put(4.1,16.1){$D$}
\put(4.0,13.3){$E$}
\put(0.5,13.3){$B$}

\put(16.0,16.0){\line(-1,1){1}}
\put(16.0,16.0){\line(-1,-1){1}}
\put(18.0,16.0){\line(1,1){1}}
\put(18.0,16.0){\line(1,-1){1}}
\put(17.0,16.0){\line(1,0){1}}
\put(17.0,16.0){\line(-1,0){1}}
\put(17.0,16.0){\line(0,1){1}}
\put(14.5,17.1){$A$}
\put(16.8,17.1){$C$}
\put(19.1,17.1){$D$}
\put(19.0,14.3){$E$}
\put(14.5,14.3){$B$}

\put(2.0,4.0){\line(-1,1){1}}
\put(2.0,4.0){\line(1,1){1}}
\put(2.0,2.0){\line(-1,-1){1}}
\put(2.0,2.0){\line(1,-1){1}}
\put(2.0,3.0){\line(1,0){1}}
\put(2.0,3.0){\line(0,1){1}}
\put(2.0,3.0){\line(0,-1){1}}
\put(0.5,5.2){$A$}
\put(3.0,5.2){$C$}
\put(3.1,2.8){$D$}
\put(2.9,0.2){$E$}
\put(0.5,0.2){$B$}

\put(12.0,15.0){\circle*{0.2}}
\put(4.0,7.0){\circle*{0.2}}
\put(7.0,4.0){\circle*{0.2}}
\put(15.0,12.0){\circle*{0.2}}
\put(12.2,15.2){$x_1$}
\put(3.8,6.3){$y_1$}
\put(6.5,3.3){$y_2$}
\put(15.2,12.2){$x_2$}

\put(4.0,7.0){\line(1,1){8}}
\put(10.0,9.0){\line(5,3){5}}
\put(10.0,9.0){\line(-3,-5){3}}

\end{picture}
\end{center}
\caption{Geodesics in tree-space consist of line segments through different
topological quadrants.
On five taxa there are 15 different quadrants, but only three are shown
above, each with a sketch of the corresponding topology.
Each axis corresponds to the length of a different split.
The shaded region does not correspond to a valid quadrant since the
splits AB$|$CDE and BE$|$ACD are incompatible.
The geodesic between $x_1$ and $y_1$ passes through three quadrants,
whereas the geodesic between $x_2$ and $y_2$ passes through just two quadrants.
In this case the geodesic is the same as the cone path.}
\label{figgeodesicandcone}
\end{figure}
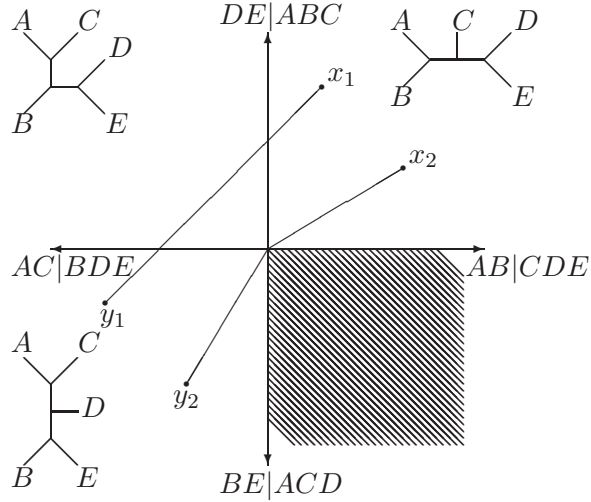

Geodesics in $\treespace$ have the following properties.
First, if $x,y\in\treespace$ have the same topology $t$, then the
geodesic joining them is the obvious Euclidean line segment in~$Q_t$.
Second, when $x$ and $y$ have some but not all splits in common, the
splits in the intersection $T_x\cap T_y$ are all included at every
point along the geodesic.
The length of the branch associated to $p\in T_x\cap T_y$ changes in
the obvious linear way from $\lambda_x(p)$ to $\lambda_y(p)$.
Third, when $x$ and~$y$ have different topologies, the geodesic may
pass through other topological orthants than the two associated with
$x$ and~$y$, as illustrated by Figure \ref{figgeodesicandcone}.
This is the case for points $x_1$ and $y_1$ in the figure.
\mbox{However}, trees along the geodesic only ever contain splits from
$T_x\cup T_y$, albeit in different combinations.
It follows that when $x$ and $y$ have different topologies, computing
the geodesic distance~$d(x,y)$ is nontrivial.
However, an efficient polynomial-time algorithm has been developed for
constructing geodesics \cite{owen09}, and we use this algorithm to
calculate distances in $\Phi$PCA.

A crucial feature of CAT$(0)$ spaces is that paths which are everywhere
locally geodesic are necessarily globally geodesic (see \cite{owen09},
Lemma 2.1). Geodesics like that between $x_1$ and $y_1$ in
Figure \ref{figgeodesicandcone} must therefore not ``bend'' as they
cross between\vadjust{\goodbreak} different orthants. For some pairs $x,y$, however, the
shortest path is given by collapsing branch lengths for splits in
$T_x\setminus T_y$ down to zero, so that the topology is then $T_x\cap
T_y$ followed by expanding out branch lengths in $T_y\setminus T_x$ to
obtain $y$. Any two points can be joined by such a path, and they are
referred to as \textit{cone paths}. In Figure \ref{figgeodesicandcone}
the cone path coincides with geodesic for points $x_2$ and $y_2$;
geodesic between $x_1, y_1$ is clearly shorter than the cone path.


\vspace*{-2pt}\section{Methods}\label{secmethods}\vspace*{-2pt}

\subsection{Existence of principal paths}

We now have the geometrical ingredients needed to define the $\Phi$PCA
procedure.
$\Phi$PCA seeks to construct a~principal path from the set of
$\treespace
$-lines defined as follows.\vspace*{-2pt}
\begin{defn}\label{defline}
A path $\Gamma$ in $\treespace$ is a $\treespace$-\textit{line} if:
\begin{longlist}
\item every sub-path of $\Gamma$ is the geodesic between its
endpoints, and
\item$\Gamma$ extends to infinity in two directions.\vspace*{-2pt}
\end{longlist}
\end{defn}

We will often just use the term \textit{line} to mean a $\treespace$-line
where the context is obvious.
Results in \cite{bill01} show that any geodesic can be extended into a~%
$\treespace$-line (though often not uniquely).
The following proposition establishes existence and uniqueness of
closest points on lines.\vspace*{-2pt}
\begin{prop}\label{propproj}
Given a $\treespace$-line $L$ and a point $x\in\treespace$, there is a~%
unique closest point $y\in L$ to $x$.\vspace*{-2pt}
\end{prop}
\begin{pf}
The proof relies mainly on the CAT$(0)$ property to
enable comparison with Euclidean space.
Let $x_0$ be any point on $L$ and suppose~$L(t)$ is a~path-length
parameterization of $L$ such that $L(0)=x_0$.
Defining $r = d(x_0,x)$, consider the triangle $x_0, x, L(t)$ for some $t>r$.
The ``skinny'' triangle property implies that
\[
d(x,L(r))< d(x,L(t)).
\]
The same bound applies to $L(-t)$.
The closest point $y\in L$, if it exists, must therefore lie on $L(t)$
for $t\in[-r,r]$.
Since this is a compact set and since the geodesic distance is a
continuous function, $d(x,L(t))$ achieves its minimum on the interval.
To prove uniqueness of the closest point $y$, suppose two distinct
points $y,y'\in L$ achieve the same minimum distance $\rho$.
Again, the ``skinny'' triangle property for the triangle $y,x,y'$ implies
that points on $L$ between $y$ and $y'$ are closer to $x$ than distance
$\rho$.
This is a contradiction, so~$y$ is unique.\vspace*{-2pt}
\end{pf}

Now suppose we are given a set of points $x_0, x_1, \ldots, x_n\in
\treespace$.
For every line~$L$ through $x_0$ we can obtain the projection
$y_1,\ldots,y_n$ of $x_1, \ldots,x_n$ onto~$L$.
This defines two functions, $f_\para(L)$ and $f_\perp(L)$, which are,
respectively, defined as the sum of squared distance along $L$, $\sum
d(x_0,y_i)^2$, and the sum of squared distances perpendicular to $L$,
$\sum d(x_i,y_i)^2$.\vadjust{\goodbreak}
\begin{prop}
There is a $\treespace$-line through $x_0$ which maximizes $f_\para$.
Similarly, there is a $\treespace$-line through $x_0$ which minimizes
$f_\perp$.
\end{prop}
\begin{pf}
We know from the proof of Proposition \ref{propproj}
that given any line $L$ through $x_0$, the points $y_i$ are at most
distance $R$ from $x_0$, where $R = \max\{d(x_0,\allowbreak x_i)\}$.
Let $S_R$ be the sphere $\{z\in\treespace\dvtx d(z,x_0)=R\}$.
Each pair $(z,z')\in S_R\times S_R$ represents a pair of geodesics
$\gamma(z,x_0)$ and $\gamma(x_0,z')$.
If $d(z,z')=2R$, then necessarily the geodesic between $z$ and $z'$ is
exact\-ly~$\gamma(z,x_0)$ followed by $\gamma(x_0,z')$, and we say $z,z'$
are \textit{antipodal}.
Every line~$L$ through $x_0$ determines an antipodal pair $(z,z')$, and
since the projected points $y_i$ all lie on the geodesic between that
pair, $f_\para$ and $f_\perp$ only depend on the pair $(z,z')$.
By continuity of the function $d\dvtx S_R\times S_R\rightarrow\mathbb{R}$,
the set of antipodal pairs is a closed subset of $S_R\times S_R$ and is
therefore compact.
It follows that there is a geodesic between antipodal points on $S_R$
which optimizes either $f_\para$ or $f_\perp$.
The geodesic can be extended into a line, and that establishes the
proposition.
\end{pf}

The optimal line may be nonunique for two reasons.
First, different extensions of the geodesic between an antipodal pair
$(z,z')$ might exist.
This would arise, for example, if all the points $x_1, \ldots, x_n$ lay
in the same topological orthant.
Second, as in regular Euclidean PCA, the collection of points $x_1,
\ldots, x_n$ can be isotropic, so that the optimal pair $(z,z')$ is nonunique.

Given the existence of optimal $\treespace$-lines, we can now consider
how to construct a principal line.
As outlined in the \hyperref[intro]{Introduction},
construction of the principal line
consists of the following steps:

{\renewcommand\thelonglist{(\arabic{longlist})}
\renewcommand\labellonglist{\thelonglist}
\begin{longlist}
\item Given trees $x_1, \ldots, x_n$, construct a ``central point'' $x_0$.
\item Given a line $L$ through $x_0$, ``project'' $x_1, \ldots, x_n$ onto
$L$ by finding the closest point $y_i$ in $L$ to $x_i$ for $i=1,\ldots,n$.
\item Find the line $L$ such which optimizes the particular choice of
objective function $f$ (either $f_\para$ or $f_\perp$).
\end{longlist}}

\noindent The details of each of these steps is described in turn, but step 3
forms the main challenge.

\subsection{Centroids and consensus}

Ideally, $x_0$ should be chosen so as to minimize the sum of squared distances:
%
%
\begin{equation}\label{equmidpoint}
x_0=\mathop{\arg\min}_{x}\sum d(x,x_i)^2.
\end{equation}
In a Euclidean vector space, this reduces to finding the mean of the
data $x_1,\ldots,x_n$.
In tree-space, due to the lack of additive structure, the mean does not
make sense, and there is no known closed solution to (\ref{equmidpoint}).
Instead, Billera et al. \cite{bill01} suggest using the \textit
{centroid}, which is defined via a~recursive procedure based on finding
the midpoint along the geodesic between any two points.
However, for large data sets,\vadjust{\goodbreak} this procedure is computationally demanding.
We therefore propose taking $x_0$ to be the \textit{majority consensus}
tree \cite{bar86}.
Finding an ``average'' or consensus tree is a well-studied problem in
phylogeny \cite{bry03} and various forms of consensus tree exist.
The majority consensus topology consists of splits which are found in
strictly more than half the trees $x_1, \ldots, x_n$.
Branch lengths on $x_0$ are assigned their average value in the data set:
\[
\lambda_{x_0}(p) = \frac{1}{| I(p) |}\sum_{i\in I(p)}\lambda_{x_i}(p)
\]
for all $p\in T_{x_0}$ where $I(p)$ is the set $\{i\dvtx p\in T_{x_i}\}$.
Results obtained later in this paper were obtained using this choice of $x_0$.
However, construction of the principal line $L$ does not rely on any
particular properties of the point~$x_0$, and $\Phi$PCA works with any
sensible choice.

\subsection{``Projection'' onto $\treespace$-lines}\label{secproj}

Given any line $L$ and points $x_1,\ldots, x_n$, Proposition \ref
{propproj} established the existence of closest points $y_1,\ldots
,y_n\in L$.
Here we describe computational aspects of this ``projection'' onto $L$.
Although this relies on existing mathematics, as presented in
Section \ref{secbackground}, the details of an algorithm for projection
onto a geodesic path have not previously been given.
We will assume $L(t)$ is a path-length parameterization of $L$.
For each point $x_i$, Euclidean projection under the embedding into
$\mathbb{R}^M$ described in Section \ref{secsplits} is used to obtain a
first guess $L(s_i)$ for $y$.
Amenta et al. \cite{amen07} showed that the geodesic distance between
two points is bounded by the Euclidean distance:
\[
\| \bolds\lambda_x - \bolds\lambda_y \|_2\leq d(x,y)\leq
\sqrt
{2}\times\| \bolds\lambda_x - \bolds\lambda_y \|_2 .
\]
It follows that if $\varepsilon_i$ denotes the Euclidean distance between
$x_i$ and $L(s_i)$, then $y_i$ lies on $L(s_i\pm\varepsilon_i)$.
This bounding interval for $y_i$ is used as the starting point for a
golden-ratio search, which is iterated until some tolerance on $y_i$ is
achieved.
It can be shown that finding $y_i$ is a convex optimization, so the
golden-ratio search is guaranteed to converge.
The proof of convexity relies on the CAT(0) property and convexity of
the equivalent Euclidean problem.
The algorithm of Owen and Provan \cite{owen09} is used to calculate
geodesic distance during the golden-ratio search.
However, it is not necessary to recompute geodesics from scratch at
every iteration: the sequence of orthants for a~geodesic at one
iteration can often be reused in the next iteration, with an associated
gain in computational efficiency.

In Euclidean vector spaces, Pythagoras' theorem gives a decomposition
of the total sum of squared distances $d_{0}^2=\sum d(x_0,x_i)^2$ of a
collection of points into contributions from directions perpendicular
and parallel to any given line $L$.
However, this decomposition does not apply in $\treespace$ with the
geodesic metric.
Nonetheless, we can evaluate the quantities
\[
d^2_\perp= \sum d(x_i,y_i)^2 \quad\mbox{and}\quad
d^2_\para= \sum d(y_i,x_0)^2
\]
for any metric.
It can be shown that for the geodesic metric, unlike the Euclidean
case, the sum of these two quantities depends on $L$.
Despite this, when evaluated for a principal path $L$, the sums of
squared distances provide a useful quantification of variability, as we
demonstrate in the results sections.

\subsection{Lines through $x_0$}

We need to construct $\treespace$-lines through $x_0$ and identify one
which optimizes our choice of objective function, $f$.
This is a~challenging problem which has not previously been considered
in the literature.
In order to achieve computational tractability, we restrict to a
particular class of $\treespace$-lines and then employ different
optimization algorithms to search over the restricted class.
To motivate this approach, we start by considering properties of lines
through $x_0$.

In the topological orthant containing $x_0$, any $\treespace$-line $L$
consists of a~straight line segment.
For brevity, we will write $T_0$ for the midpoint topology~$T_{x_0}$
and~$\lambda_0$ for the branch length function $\lambda_{x_0}$.
Let~$Q_0$ denote the orthant containing $x_0$, and for now assume that~%
$x_0$ is a binary tree (so it contains the maximal number of splits).
If $p$ is a split contained in $Q_0$, then the branch length at a point
$y(s)\in Q_0$ on $L$ has the form
%
%
\begin{equation}\label{equgeocoord}
\lambda_{y(s)}(p) = \lambda_0(p) +s\times w_p,
\end{equation}
where $w_p$ is a ``weight'' associated to split $p$, and $s$ lies on some
interval containing zero.
The set of weights determines the direction vector of the line segment
through~$x_0$.
Given such a line segment, we need to know how it might extend beyond
$Q_0$ into the rest of tree-space.

Where the segment meets a face of $Q_0$, at least one split is assigned
zero branch length.
Generically, the line segment will meet a co-dimension $1$ face of
$Q_0$, so that just one split $p$ will have zero length.
Solving equation (\ref{equgeocoord}) for this split gives
%
%
\begin{equation}\label{equsolns}
\lambda_{y(s)}(p) =0 \quad\Rightarrow\quad s=-\lambda_0(p)/w_p.
\end{equation}
The line then extends from this point into one of the neighboring
alternative orthants.
In a similar way, every other split whose length varies in the initial
line segment containing $x_0$ is associated with a solution of
equation (\ref{equsolns}) and, correspondingly, with an alternative
split related to the first by NNI.
If we restrict to the set of ``generic'' lines (those which always meet a
co-dimension~$1$ face of every orthant), then finding the optimal line
$L$ therefore consists of a topological problem (namely, choosing a new
split $p'$ to replace each $p$) and a geometrical problem (finding the
best set of weights $w_p$).
However, these problems are not independent.
We can order the solutions to (\ref{equsolns}) as we move out from
$x_0$ in a particular direction along $L$.
Suppose the first solution we come to is at $s=s_1$ and we replace
split $p_1$ with $p'_1$.
At the next solution $s=s_2$, split $p_2$ is assigned zero length and
we replace it via an NNI move.
However, the choice of splits available as replacements for $p_2$
does
not depend solely on $p_2$ but also on the rest of the tree topology
just\vadjust{\goodbreak} before $s=s_2$---and therefore potentially on the choice of
replacement $p'_1$ of $p_1$.
Thus, the topological aspect of construction depends on the relative
order of the solutions to (\ref{equsolns}), which in turn
depends on the weights $w_p$.
Optimization over the set of possible weights and splits will be
computationally demanding for trees with more than a few species---an
exhaustive search will not be possible.

A key feature of the description above is the assumption that line
segments meet the boundary of orthants in co-dimension $1$ faces.
We restrict our search space for $L$ similarly, but take into account
the possibility that $x_0$ might not be fully resolved.
We make this more formal as follows.
\begin{defn}
Suppose $p\in S$ is compatible with $x_0$ and $p'\in S$ is
obtained by extended nearest neighbor interchange of $p$ in $T_0\cup\{
p\}$.
The \textit{simple line} through $x_0$ associated with $p,p'$ and weight
$w$ is the path $y(s)\in\treespace$ defined by
\begin{eqnarray*}
\lambda_{y(s)}(p) &=&
\lambda_0(p) +sw \qquad\mbox{if } \lambda_0(p) +sw \geq0\\
&=& 0 \qquad\mbox{otherwise},\nonumber\\
\lambda_{y(s)}(p') &=&
-\bigl(\lambda_0(p) +sw \bigr) \qquad\mbox{if }\lambda_0(p) +sw \leq0\\
&=&0 \qquad\mbox{otherwise},\nonumber\\
\lambda_{y(s)}(q) &=& \lambda_0(q) \qquad\mbox{if }q\notin\{p,p'\}.
\end{eqnarray*}
\end{defn}

Such a path moves through a single pair of orthants.
The next definition extends simple lines to pass through more than two orthants.
\begin{defn}\label{defext}
Suppose $x(s)$ is the simple line through $x_0$ defined by split pairs
$(p_1,p'_1),\ldots,(p_k,p'_k)$ and weights $w_1,\ldots,w_k$, and
suppose that the pair of splits $(p_{k+1}, p'_{k+1})$ and weight
$w_{k+1}\in\mathbb{R}$ satisfy the following:
\begin{longlist}
\item$p_{k+1}$ is compatible with $x(s)$ for all $s$ such that
$\lambda
_0(p_{k+1}) +sw_{k+1} \geq0$,
\item$p'_{k+1}$ is compatible with $x(s)$ for all $s$ such that
$\lambda_0(p_{k+1}) +sw_{k+1} < 0$,
\item$p_{k+1},p'_{k+1}$ are related by XNNI in $x(s_{k+1})$ where
$s_{k+1}=-\lambda_0(p_{k+1})/\allowbreak w_{k+1}$.
\end{longlist}
Then the simple line $y(s)$ defined by $(p_1,p'_1),\ldots
,(p_{k+1},p'_{k+1})$ and weights $w_1,\ldots,w_{k+1}$ is given by
%
%
\begin{eqnarray}
\label{equbl1}
\lambda_{y(s)}(p_i) &=&
\lambda_0(p_i) +s w_i \qquad\mbox{if }\lambda_0(p_i) +s w_i\geq
0\nonumber\\[-8pt]\\[-8pt]
&=&0 \qquad \mbox{otherwise},\nonumber\\
\label{equbl2}
\lambda_{y(s)}(p'_i) &=&
-\bigl(\lambda_0(p_i) +s w_i\bigr) \qquad\mbox{if }\lambda_0(p_i) +s w_i\leq
0\nonumber\\[-8pt]\\[-8pt]
&=&0 \qquad \mbox{otherwise},\nonumber\\
\label{equbl3}
\lambda_{y(s)}(q) &=& \lambda_0(q) \qquad\mbox{if }q\notin\{p_i\}\cup
\{p'_i\},\qquad i=1,\ldots,k+1.
\end{eqnarray}
\end{defn}

To prove that simple lines satisfy the conditions of Definition \ref
{defline}, Proposition 4.2 of \cite{bill01} can be applied to any pair
of points on a simple line in order to show that the subpath between
those points is geodesic.

Simple lines through $x_0$ resemble the geodesic between $x_1$ and
$y_1$ in Figure~\ref{figgeodesicandcone}: they continue between
orthants without bends, and hence are always locally geodesic.
Moreover, any straight line segment through $x_0$ can be obtained as
part of a simple line.
Nonetheless, restriction to the class of simple lines removes many
lines from consideration.
Cone paths are ruled out, together with any geodesic for which some
subset of the splits changes like a~cone path.
(This latter class of geodesics resembles $x_1-y_1$ in Figure \ref
{figgeodesicandcone} for some splits and $x_2-y_2$ for others.)
This restriction is carried out for the sake of computational tractability.
More discussion is given in Section \ref{secconc}.

Definition \ref{defext} describes how to extend a simple line on $k$
split pairs to one on $k+1$ split pairs.
Our algorithms for finding an optimal simple line are based precisely
on this operation.
Suppose a simple line $L$ is determined by sets of splits $P=\{
p_1,\ldots,p_k\},P'=\{p'_1,\ldots,p'_k\}$ and weights $W=\{w_1,\ldots
,w_k\}$.
Conditions (i)--(iii) of Definition \ref{defext} place constraints on
any proposed splits $p,p'$ and weight $w$ which might be used to extend $L$.
The values $s_i=-\lambda_0(p_i)/w_i$ correspond to points at which $L$
crosses the boundary between orthants, and we can assume they are
ordered with $s_1\leq s_2\leq\cdots\leq s_k$.
They divide $L$ up into $k+1$ intervals $I_i=[s_i,s_{i+1}]$ for
$i=0,1,\ldots,k$ taking $s_0=-\infty$ and $s_{k+1}=\infty$, such that
the topology of $L$ is constant on each interval.
Let $t_i$ denote the tree topology on $I_i$.
Now suppose that $p$ is compatible with $T_0$ and $p'$ is a proposed
replacement for~$p$.
Suppose we also propose an interval $I_i$ on which we require the XNNI
move to be performed.
Conditions (i)--(iii) are then equivalent to the following.\looseness=1

\textit{Geometrical constraint}:
\[
s_i \leq-\frac{\lambda_0(p)}{w} \leq s_{i+1}\mbox{, so the XNNI
move occurs on interval $I_i$}.
\]

\textit{Topological constraint}:
\begin{itemize}
\item If $w<0$, then $p$ must be compatible with $t_0, \ldots, t_i$ and
$p'$ must be compatible with $t_i\setminus p, \ldots, t_{k+1}\setminus
p$; or
\item if $w>0$, then $p'$ must be compatible with $t_0\setminus p,
\ldots, t_i\setminus p$ and $p$ must be compatible with $t_i, \ldots,
t_{k+1}$.
\end{itemize}

When $p$ is not contained in $T_0$, but instead \textit{extends} the
midpoint topology, then $\lambda_0(p)=0$ and the solution to
equation (\ref{equsolns}) is $s=0$.
In this case, the geometric constraint corresponds to an unbounded
interval for $w$, and the interval $I_i$ on which the XNNI move
$p\rightarrow p'$ takes place must necessarily contain $s=0$.

\subsection{Greedy algorithm for finding an optimal simple line}

The following algorithm repeatedly extends a simple line by adding in a\vadjust{\goodbreak}
new split pair at each iteration.
The pair chosen is the one which gives the best improvement in the
objective $f$:

{\renewcommand\thelonglist{(\arabic{longlist})}
\renewcommand\labellonglist{\thelonglist}
\begin{longlist}
\item Let $F$ be the set of feasible splits (see below).
\item Consider in turn every split $p$ in $F$ that is compatible with
$T_0$, and every possible replacement $p'$ for $p$.
\item For each interval $I_i$, test whether $p'$ is an XNNI replacement
of $p$ in~$t_i$.
\item If (ii) holds on interval $I_i$, then next check whether the pair
$p,p'$ satisfies either topological constraint for that interval. Fix
the sign of $w$ depending on which constraint applies.
\item If either topological constraint holds, then find $w$ that
maximizes the variance of the projected points on $L$, subject to the
geometrical constraint and sign of~$w$.
This is carried out using the golden ratio search for the optimum value
of $w$.
\item Repeat for all feasible pairs $p,p'\in F$ and find the pair that
gives the maximum projected variance.
\item Add $p$, $p'$ and $w$ to the lists $P$, $P'$ and $W$, and reorder
the lists according to the solutions of (\ref{equsolns}).
Remove $p,p'$ from $F$, and repeat from step 2.
\end{longlist}}

\noindent The algorithm continues until no more splits can be added to $L$.
This will occur in at most $m-3$ iterations where $m$ is the number of
species, since every tree can contain at most $m-3$ nontrivial splits.

The set of feasible splits $F$ could be taken to be the entire set of
possible splits~$S$, but this is inefficient.
If neither split $p,p'$ is contained in any of the trees $x_1,\ldots,
x_n$, then adding the pair to $L$ will only increase the distances
$d(x_i,y_i)$ so that $L$ is a worse approximation to the data.
We therefore take~$F$ to be the set of nontrivial splits found in at
least one tree $x_1,\ldots,x_n$.
It is possible that at some stage the best improvement in $f$ might be
given by some $p\in F$ and a replacement $p'\notin F$
(e.g., consider the case that all the trees $x_i$ lie in the
same orthant).
However, in such a situation, the data would not be informative about
the choice of $p'$, and so we disregard this possibility.

The greedy algorithm terminates after at most $m-3$ iterations.
During each iteration $O(|F|^2)$ pairs of splits and $O(m)$ possible
intervals for the move $p\rightarrow p'$ are considered.
For each pair of splits and interval, $n$ trees are projected onto the
proposed line.
Each projection requires $O(m^4)$ steps.
The golden ratio search during projection is performed to a fixed
tolerance, and so is independent of $m$, $n$ and $|F|$.
Overall, the algorithm therefore requires $O(m^6\times n\times|F|^2)$
steps where $F$ is at worst $O(nm)$.

\subsection{Monte Carlo optimization algorithm}

A Monte Carlo optimization algorithm was also implemented in order to
provide comparisons with the greedy approach.
A simulated-annealing type approach was adopted, where the state at
each iteration comprised a simple line $L$ through $x_0$.
At each iteration a ``birth'' or ``death'' move was randomly proposed from
the current state.
Birth moves consisted of adding a valid split pair to $L$, while death
moves consisted of removing a split pair from $L$.
Birth moves were obtained by selecting $p\in F$ uniformly at random,
then selecting $p'$ uniformly at random from the possible XNNI
replacements of $p$ satisfying the constraints defined above.
The weight assigned to $p,p'$ was obtained by the golden ratio search,
as for the greedy approach.
Death moves were carried out by choosing at random the split pair at
either end of $L$ (i.e., with largest positive or negative $s_i$) and
removing it.
Removing other split pairs results in incompatible sets of splits along
$L$ and is therefore forbidden.
The relative probabilities of birth and death depended on the number
$k$ of split pairs in~$L$ and were designed to favor birth for small
$k$ and death when $k$ was large.
Proposals leading to improvement in the objective were always accepted.
Other proposals were accepted with probability
\[
\operatorname{Pr}(\mbox{accept}) = \biggl(1-\frac{\delta}{D} \biggr)
^{{1}/{\tau}},
\]
where $\delta$ is the absolute difference of the objective for the
proposed and current state, $D$ is a bound for $\delta$, and $\tau$ is
the ``temperature.''
For $f=f_\perp$, $D$ was taken to be $\sum d(x_0,x_i)^2$, while for
$f=f_\para$, $D$ was taken to be $f_\para$ for the current state.
The temperature $\tau$ was slowly decreased as the optimization progressed.

\subsection{Branch length transformations}\label{sectrans}

We investigated certain transformations of the data $x_1, \ldots, x_n$
prior to analysis with $\Phi$PCA.
Kupczok et al. \cite{kup08} suggest scaling each tree $x_1, \ldots,
x_n$ to have the same total branch length.
In practice, this seemed to make little difference to the examples we
looked at in the results section below.
Instead we considered the following branch length normalization.
For each split $p$, branch lengths were scaled by a constant so that
the average branch length associated with $p$ across the whole data set
was unity.
This was repeated for each split in the data set.
The idea behind this is to make $\Phi$PCA measure variability relative to
branch length and to amplify the variability in short branches.
In regular PCA the correlation matrix can be analyzed instead of the
covariance matrix, and this branch length transformation can be thought
of as being analogous to the correlation matrix version.
Principal geodesics obtained for branch-length normalized data can be
back-transformed onto the original scale by scaling the weights $W$.


\section{Relationship to the work by Wang and Marron}

Wang and Marron~\cite{wang07} previously developed PCA in a space of
trees, and on account of the similarities of our approach to theirs, in
this section we look in\vadjust{\goodbreak} detail at the relationship between the approaches.
The steps underlying our approach specified at the start of
Section~\ref
{secmethods} were taken directly from \cite{wang07}, but the details of
how these steps are carried out are quite different on account of the
different geometries under consideration.

In \cite{wang07} rooted bifurcating trees are considered, but, unlike
phylogenetic trees, the leaf vertices are not assigned taxon labels.
Instead, each vertex can have a ``left'' and a ``right'' descendant, and
trees in the data set can have different depths from root to leaf.
Most importantly, branches do not have any associated length, but,
instead, each vertex present in a tree has an associated real number
(or vector).
An example of such data consists of blood vessel information from
medical imaging: vertices represent blood vessels, edges represent
connections between blood vessels, and the data value associated to
each vertex corresponds to some measurement at that point in the blood
vessel structure.
One crucial difference between the two spaces of trees is that in \cite
{wang07} there is no relationship between the values associated to
vertices and the topological structure of the tree.
This is different from the space $\treespace$, in which branch lengths
can be shrunk down and replaced by an alternative topology.

This separation of ``topological'' and ``geometrical'' aspects of the
problem in \cite{wang07} results in principal components with separate
topological and geometrical parts.
In Wang and Marron's terminology, a \textit{structure tree line} is a
sequence of vertices, each descended from the previous vertex, which
can be thought of as (discontinuous) ``growth'' of a tree toward a leaf,
by grafting on branches.
In contrast, an \textit{attribute tree line} consists of a fixed tree
structure with ``direction vectors'' associated to vertices.
This is clearly very different from the lines constructed by $\Phi$PCA in
which the principal path reflects both topological and geometrical
variability in the data set.

Not surprisingly, given the different structures of the spaces
considered, the metrics used in the two approaches differ.
The metric in \cite{wang07} is a linear combination of the (unweighted)
Robinson--Foulds metric \cite{rob81} and a Euclidean distance between
the vectors associated to each vertex.
This metric is inexpensive to compute, in contrast to the geodesic
metric which we consider, and this reduces the computational burden of
their approach relative to ours.
The midpoint $x_0$ in \cite{wang07} is taken to have the majority
consensus topology \cite{bar86}, as used here, since this minimizes the
sum of the Robinson--Foulds distances of the midpoint from $x_1,\ldots,x_n$.
However, while Wang and Marron obtain an exact form of Pythagoras'
theorem with their metric, that is not the case for $\Phi$PCA (see
Section \ref{secproj}).

In summary, the method of Wang and Marron cannot be applied directly to
phylogenetic trees, since the trees they consider lack taxon labels and
branch lengths and so cannot represent phylogenies.
While our approach builds on the same framework as that laid out in
\cite{wang07}, differences in the geometries of the spaces under
consideration make the mathematical details of the implementation of
PCA substantially different.
In particular, $\Phi$PCA relies heavily on the geometry of $\treespace$
described by Billera et al. \cite{bill01}.
It is interesting to note how a seemingly small difference in the
geometry of the space under consideration can substantially change the
way PCA is implemented.


\section{Simulation studies}\label{secsims}

\subsection{Simple mixtures}\label{sec41}

$\Phi$PCA was used to analyze collections of randomly generated trees
with two (or more) known underlying topologies.
These simulations were not intended as a model of a specific process
giving rise to alternative trees, but were performed in order to verify
the methodology and demonstrate how it works on simple examples.
We describe the simulations very briefly here, but give more details in
the supplementary material~\cite{Nye11}.
Two sets of simulations were performed.
In the first set, trees were simulated such that each had one of two
possible topologies $t_1$ or $t_2$.
The underlying topologies $t_1,t_2$ were related by an NNI move and
represented alternative positions for a clade within the tree.
Topology $t_1$ was adopted with probability~$\theta$ and $t_2$ with
probability $1-\theta$.
Apart from branches affected by the change in topology, all other
branch lengths were kept fixed.
For each value of~$\theta$, 100 trees were randomly generated in this way.
Additional variability was added by simulating a DNA alignment for each
tree, and then replacing the tree with the maximum likelihood (ML) tree
estimated from the alignment.
$\Phi$PCA was used to analyze these estimated trees.
A~second set of simulations was performed in which there were two
correlated changes in topology.
Each tree consisted of two subtrees, and each subtree had either
topology~$t_1$ or $t_2$ as in the first set of simulations.
The alternative topologies in each half of the tree were simulated to
arise with correlation~$\rho$.
Again, additional variability was added by simulating alignments and
replacing each tree with an ML estimate.
100 trees were generated for each pair of values~$\theta, \rho$ and
$\Phi$PCA was used to analyze each set of estimated ML trees.

The results indicated that optimization of $f_\para$ gave the best
performance: paths obtained by optimizing $f_\perp$ sometimes missed
the changes in topology imposed in the data sets.
In this non-Euclidean setting the sum $d^2_\para+d^2_\perp$ is
generally less than the total sum of squared distances $d^2_0$, so
optimization of~$f_\perp$ may result in principal paths that fail to
capture variability in the data by finding paths in which both sums
$d^2_\para$ and $d^2_\perp$ are small.

In both sets of simulations, $\Phi$PCA with $f_\para$ gave principal
paths corresponding to the change between the imposed alternative topologies.
In the first set of simulations, based on a single pair of alternative
topologies, as~$\theta$ increased the change between the underlying
topologies $t_1,t_2$ dominated the principal path (the corresponding
splits received a higher weight) as variability due to tree estimation
from alignments was dominated by the imposed variability in
topology.\vadjust{\goodbreak}
In the second set of simulations, for small values of $\rho$ the
principal path corresponded to change between the alternative
topologies on one part of the tree, with the other pair of alternative
topologies receiving a low weight.
For larger $\rho$, the correlated alternatives in both parts of the
tree were identified by the principal path.
More details are given in the supplementary material
\cite{Nye11}.

\subsection{Long branch attraction}

\begin{figure}

\includegraphics{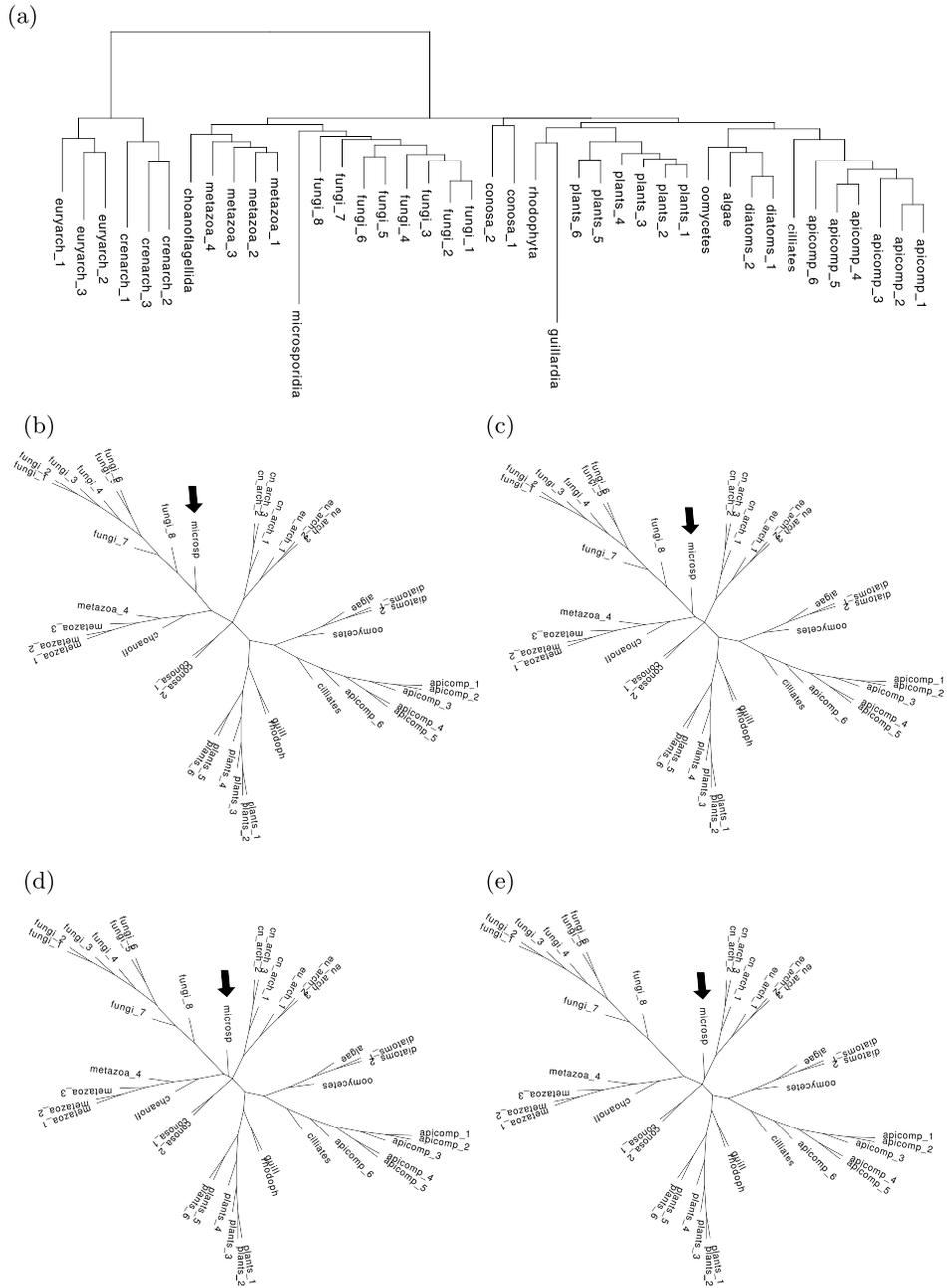}

\caption{Simulation study of LBA.
\textup{(a)} Underlying tree with two long branches and distant out-group
(archaea, on the left).
\textup{(b)--(e)} Trees along the principal path.
Branches were normalized to have unit mean.
No back-transform to the original scale was performed, since this
obscured the visual effect.
Arrows highlight the microsporidia group (labeled ``microsp'') moving
round to join the outgroup (labeled ``cn\_arch'' and ``eu\_arch'').}
\label{figLBA}
\end{figure}

In order to demonstrate a potential application of $\Phi$PCA, a simple
study of long branch attraction (LBA) was performed.
LBA is a feature of phylogenetic methods in which species on long
branches are often grouped together erroneously on estimated trees.
We took a tree from the literature \cite{bri05} representing a deep
phylogeny of eukaryote species which includes two long branches and a
distant out-group, as shown in Figure \ref{figLBA}.
100 trees were simulated by first simulating amino acid alignments from
the base tree (300 base pairs, $\mbox{WAG}\mbox{$+$}4\Gamma$ model) using the seq-gen
software \cite{ram97} and then obtaining an ML estimate tree from each
alignment using phyML \cite{guin03}.
$\Phi$PCA was used to analyze the set of trees estimated from the
simulated alignments.

%
%

Analysis of the simulated trees was carried out first with
un-normalized data and then again with the normalization procedure
described in Section \ref{sectrans}. Optimization was carried out using
the $f_\para$ objective function, and results were obtained with both
the greedy and Monte Carlo algorithms. Despite long runs, the Monte
Carlo algorithm failed to improve on the results obtained with the
greedy algorithm. The results obtained with the greedy algorithm are
shown in Table \ref{tabLBA} and Figure \ref{figLBA} shows the
principal path obtained with normalized branch lengths. The
``proportion of variance'' $d^2_\para/d^2_0$ was greater\vspace*{1pt}
with the normalized data, and so we suggest that normalization is
preferable for this data set.

%
%
\begin{table}
\caption{Results of the LBA simulations. The largest two weights $w$
for split pairs and a~description of the corresponding changes in
topology are given. Weights were normalized to have unit Euclidean
norm}\label{tabLBA}
\begin{tabular*}{\tablewidth}{@{\extracolsep{\fill}}ll@{}}
\hline
$\bolds{w}$ & \multicolumn{1}{c@{}}{\textbf{Change in topology}} \\
\hline
\multicolumn{2}{@{}c@{}}{\textit{Without normalization}
($d^2_\para/d^2_0=2.6\%$)} \\[2pt]
$0.859$ & Guilardia moves past pairing with Rhodophyta to top of clade
with plants\\
$0.152$ & Guilardia moves from top of clade with plants to position
closer to Archaea
\\[4pt]
\multicolumn{2}{@{}c@{}}{\textit{Branch lengths normalized}
($d^2_\para/d^2_0=10.3\%$)}
\\[4pt]
$0.706$ & Microsporidia moves from grouping with fungi to top of clade
with Metazoa\\
$0.482$ & Microsporidia grouped with Archaea\\
\hline
\end{tabular*}
\end{table}

As explained in
\cite{bri05}, estimated trees tend to place the long
branches (\textit{Guillardia} and Microsporidia) next to the out-group
(Archaea).
The analyses of both the\vadjust{\goodbreak} un-normalized and normalized data show this
effect with, respectively, \textit{Guillardia} and Microsporidia
``floating'' round the tree to be placed closer to the Archaea.
The fact that each principal path captures a single such effect
suggests that the attraction of the two long branches to the Archaea is
uncorrelated in the data.
$\Phi$PCA exactly captures the expected LBA artefact in the simulated data.


\section{Analysis of metazoan data}\label{secresults}

$\Phi$PCA was applied to a set of 118 gene trees from 21 metazoan
(animal) species, previously analyzed in \cite{kup08}.
$\Phi$PCA was performed on both unscaled and branch length normalized
data using the $f_\para$ objective function.
The Monte Carlo optimization algorithm obtained principal paths with
slightly higher $f_\para$ scores than the greedy algorithm, and so we
refer to that set of results here.
The principal paths obtained with the two algorithms were similar, and
shared the majority of split pairs in common.
The ``proportion of variance'' $d^2_\para/d^2_0$ was $1.8\%$ for the
unscaled data and $4.6\%$ for the normalized data---relatively low in
both cases.
However, the simulation studies produced similarly low scores (between
$3\%$ and $5\%$ on artificial data), suggesting that low scores might
be common even when $\Phi$PCA is successfully capturing aspects of
variability in the data.
Further comments about the low proportion of variance are made in
Section \ref{secconc}.

\begin{figure}

\includegraphics{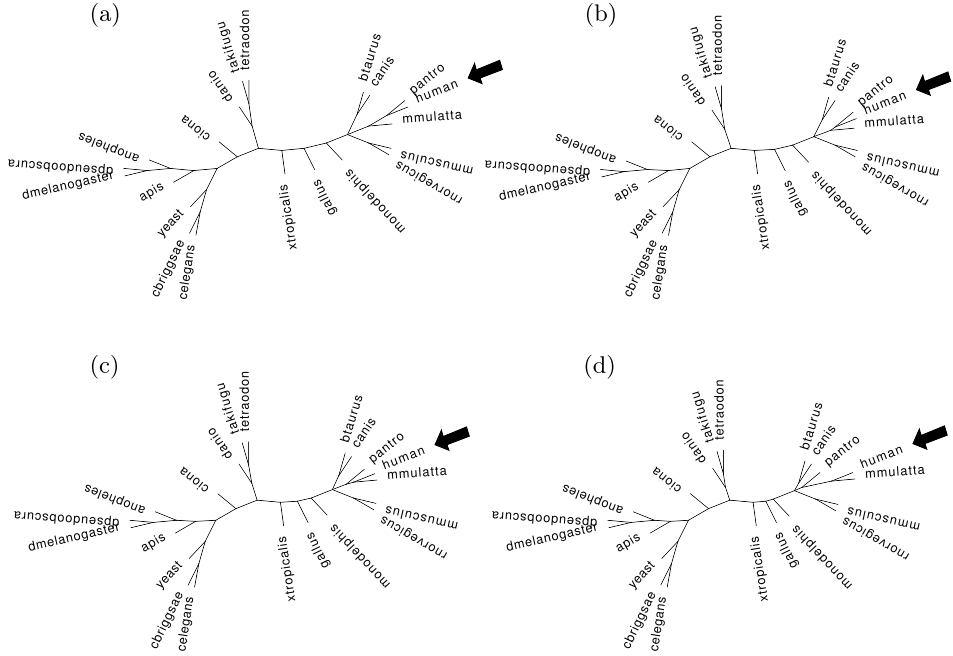}

\caption{Trees along the principal path for the normalized metazoan
data in order \textup{(a)--(d)}.
\textup{(a)}~corresponds to the majority consensus topology.
Human moves from being grouped with Chimp (labeled ``pantro'') to Macaque
(labeled ``mmulatta''), as highlighted by the solid arrows.}
\label{figmetazoa}
\vspace*{-4pt}
\end{figure}

The principal path obtained for the unscaled data corresponded to
uncertainty in the positioning of the out-group, yeast.
It moves from being placed next to the worms to being grouped with sea squirt.
This might be an LBA effect since sea squirt and yeast lie on
relatively long branches.
Results of the analysis using data with normalized branch lengths are
shown in Figure \ref{figmetazoa}.
%
%
The principal path indicates uncertainty in the placement of Human: it
is either grouped with Chimpanzee or Macaque.
The position of Human relative to its neighbors was a longstanding
problem in phylogenetics~\cite{good98}.
Uncertainty in the positioning arises from the presence of a relatively
short internal branch in the species tree joining Human and Chimpanzee
to the other primates.
Although well known in evolutionary biology, this simple example
illustrates how $\Phi$PCA can be used to identify and visualize
alternatives within a set of trees.


\section{Conclusion}\label{secconc}

We have presented a procedure for identifying principal paths in the
space of phylogenetic trees which best approximate a set of alternative
phylogenies in an analogous way to standard PCA in Euclidean vector spaces.
A key feature of the approach is the use of metrics that combine
geometric and topological information about trees.
The principal paths constructed coincide with the geodesic between
every pair of points on the path.
Each principal path is equipped with a summary statistic analogous to
the Euclidean proportion of variance which quantifies variability along
the path.

Results obtained from simulated and experimental data sets gave
values
for the ``proportion of\vadjust{\goodbreak} variance'' $d_\para^2/d_0^2$ which were
relatively low in comparison to typical values for standard PCA (e.g.,
about $5\%$ for the normalized metazoa data set).
This is the result of two features of the problem.
First, the data sets analyzed in this paper are high dimensional
(containing over $100$ different splits), and in the same way as for
standard Euclidean data, this tends to lead to lower proportions of variance.
To illustrate this, consider, for example, a multivariate normal
distribution with dimension $100$ and covariance matrix $
\operatorname{diag}(5\sigma^2, \sigma^2,\ldots, \sigma^2)$.
Standard PCA would give a proportion of variance of roughly $5\%$ even
though the variance along the principal component is substantially
higher than in other directions.
Second, the failure of Pythagoras' theorem in tree-space means that for
any analysis variance ``leaks out,'' that is, the sum of squares
$d_\parallel^2+d_\perp^2$ is less than the total sum of square
distances for the original data $d_0^2$, further decreasing $d_\para^2/d_0^2$.

In order to construct principal paths, two approximations have been imposed:

{\renewcommand\thelonglist{(\arabic{longlist})}
\renewcommand\labellonglist{\thelonglist}
\begin{longlist}
\item A greedy algorithm or Monte Carlo search is carried out in the
configuration space of paths in order to find the optimal path.
There is no guarantee that the optimal path will always be
found.\vadjust{\goodbreak}
\item The configuration space itself is restricted to a subclass of
paths (referred to as simple lines).
By restricting in this way we might rule out capturing types of
variability in the data under analysis.
\end{longlist}}

\noindent
We consider the second of these approximations to be more limiting, and
it is difficult to generalize the approach we have described to
overcome it.
Another area requiring further research is the construction of
higher-dimensional approximations to data, analogous to second and
third, etc. principal components.
Construction of higher order principal components in Euclidean vector
spaces is carried out by working orthogonally to the first principal component.
Tree-space is not equipped with an inner product, and so this procedure
cannot carry over directly to~$\treespace$.
Our algorithms for constructing the principal path $L$ do not therefore
readily generalize to give higher order paths.
Instead, we would need to consider two-dimensional subsets of
$\treespace$ which approximate the data $x_1, \ldots, x_n$ as closely
as possible.
In analogy to the definition of $\treespace$-lines, we would require
that such subsets~$\Pi$ contain the geodesic between any two points in
$\Pi$, and so $\Pi$ would locally resemble a plane in each orthant.
However, in contrast to the theory of geodesics, the theory of
higher-dimensional surfaces in tree-space is not well developed.
We have not attempted to advance this theory in this paper, but have
focused on the already considerable problem of identifying lines which
best approximate the data.

The $\Phi$PCA procedure has been presented as an empirical analysis of
sampled trees without reference to any underlying distribution that
generated the trees.
Distributions such as sampling distributions, bootstrap distributions
and Bayesian posteriors are of fundamental importance in phylogenetic
inference, but the geometrical properties of these distributions have
received little study.
Billera et al. \cite{bill01} considered spherically symmetric
distributions with density decaying exponentially away from a central point.
A second form of isotropic distribution consists of the limit of a
random walk in tree-space from a fixed central point.
By simulating samples from a~suitable random walk and carrying out
$\Phi$PCA on the samples, an empirical $p$-value could be assigned to the
proportion of variance of a principal line constructed from
experimental data, as a test for significant departure from isotropy.
Holmes \cite{hol05}, however, suggests that the assumption of spherical
symmetry is not realistic for most distributions of interest.
One area where distributions on tree-space have been defined more
precisely is the study of population-genetic effects on gene
phylogenies \cite{deg05}.
Such distributions could be studied in the context of tree-space
geometry, and it might be possible to obtain the sampling theory of
principal lines under $\Phi$PCA in this case.

This paper has presented the results of applying $\Phi$PCA to some
relatively simple examples, and demonstrated the type of information
principal paths reveal.
The method can be applied to larger data sets and it has the potential
to provide new insights into a range of problems in evolutionary biology.
Software for performing $\Phi$PCA and for visualizing principal paths as
animations of trees is available in the supplementary material \cite{Nye11},
together with the data sets analyzed in this paper.


\section*{Acknowledgments}

The author would like to thank Anne Kupczok and the group of Arndt von
Haeseler for generously providing the set of metazoan gene trees. I am
also grateful to two anonymous reviewers for their very helpful
comments.

\begin{supplement}[id=suppA]
\stitle{Principal components analysis in the space
of phylogenetic trees: Supplementary information}
\slink[doi]{10.1214/11-AOS915SUPP} 
\sdatatype{.pdf}
\sfilename{aos915\_supp.pdf}
\sdescription{This contains further information about the
simulation studies in Section \ref{sec41}.}
\end{supplement}

%

\printaddresses


\begin{thebibliography}{29}

\bibitem{all01}
%
\begin{barticle}[mr]
\bauthor{\bsnm{Allen},~\bfnm{Benjamin~L.}\binits{B.~L.}} \AND
\bauthor{\bsnm{Steel},~\bfnm{Mike}\binits{M.}}
(\byear{2001}).
\btitle{Subtree transfer operations and their induced metrics on evolutionary
trees}.
\bjournal{Ann. Comb.}
\bvolume{5}
\bpages{1--15}.
\bid{doi={10.1007/s00026-001-8006-8}, issn={0218-0006}, mr={1841949}}
\bptok{imsref}%
\end{barticle}
%
\endbibitem

\bibitem{amen07}
%
\begin{barticle}[mr]
\bauthor{\bsnm{Amenta},~\bfnm{Nina}\binits{N.}},
\bauthor{\bsnm{Godwin},~\bfnm{Matthew}\binits{M.}},
\bauthor{\bsnm{Postarnakevich},~\bfnm{Nicolay}\binits{N.}} \AND
\bauthor{\bsnm{St.~John},~\bfnm{Katherine}\binits{K.}}
(\byear{2007}).
\btitle{Approximating geodesic tree distance}.
\bjournal{Inform. Process. Lett.}
\bvolume{103}
\bpages{61--65}.\break
\bid{doi={10.1016/j.ipl.2007.02.008}, issn={0020-0190}, mr={2322069}}
\bptok{imsref}%
\end{barticle}
%
\endbibitem

\bibitem{ayd09}
%
\begin{barticle}[mr]
\bauthor{\bsnm{Aydin},~\bfnm{Burcu}\binits{B.}},
\bauthor{\bsnm{Pataki},~\bfnm{G{\'a}bor}\binits{G.}},
\bauthor{\bsnm{Wang},~\bfnm{Haonan}\binits{H.}},
\bauthor{\bsnm{Bullitt},~\bfnm{Elizabeth}\binits{E.}} \AND
\bauthor{\bsnm{Marron},~\bfnm{J.~S.}\binits{J.~S.}}
(\byear{2009}).
\btitle{A principal component analysis for trees}.
\bjournal{Ann. Appl. Stat.}
\bvolume{3}
\bpages{1597--1615}.
\bid{doi={10.1214/09-AOAS263}, issn={1932-6157}, mr={2752149}}
\bptok{imsref}%
\end{barticle}
%
\endbibitem

\bibitem{bar86}
%
\begin{barticle}[mr]
\bauthor{\bsnm{Barth{\'e}l{\'e}my},~\bfnm{Jean-Pierre}\binits{J.-P.}}
(\byear{1986}).
\btitle{The median procedure for {$n$}-trees}.
\bjournal{J. Classification}
\bvolume{3}
\bpages{329--334}.
\bid{doi={10.1007/BF01894194}, issn={0176-4268}, mr={0874243}}
\bptok{imsref}%
\end{barticle}
%
\endbibitem

\bibitem{bill01}
%
\begin{barticle}[mr]
\bauthor{\bsnm{Billera},~\bfnm{Louis~J.}\binits{L.~J.}},
\bauthor{\bsnm{Holmes},~\bfnm{Susan~P.}\binits{S.~P.}} \AND
\bauthor{\bsnm{Vogtmann},~\bfnm{Karen}\binits{K.}}
(\byear{2001}).
\btitle{Geometry of the space of phylogenetic trees}.
\bjournal{Adv. in Appl. Math.}
\bvolume{27}
\bpages{733--767}.
\bid{doi={10.1006/aama.2001.0759}, issn={0196-8858}, mr={1867931}}
\bptok{imsref}%
\end{barticle}
%
\endbibitem

\bibitem{bri05}
%
\begin{barticle}[author]
\bauthor{\bsnm{Brinkmann},~\bfnm{H.}\binits{H.}},
\bauthor{\bparticle{{v}an~{d}er} \bsnm{{G}eizen},~\bfnm{M.}\binits{M.}},
\bauthor{\bsnm{Zhou},~\bfnm{Y.}\binits{Y.}},
\bauthor{\bparticle{{P}oncelin~{d}e} \bsnm{{R}aucourt},~\bfnm
{G.}\binits{G.}}
\AND\bauthor{\bsnm{Philippe},~\bfnm{H.}\binits{H.}}
(\byear{2005}).
\btitle{An empirical assessment of long-branch attraction artefacts in deep
eukaryotic phylogenomics}.
\bjournal{Syst. Biol.}
\bvolume{54}
\bpages{743--757}.
\bptok{imsref}%
\end{barticle}
%
\endbibitem

\bibitem{bry03}
%
\begin{bincollection}[mr]
\bauthor{\bsnm{Bryant},~\bfnm{David}\binits{D.}}
(\byear{2003}).
\btitle{A classification of consensus methods for phylogenetics}.
In \bbooktitle{Bioconsensus ({P}iscataway, {NJ}, 2000/2001)}.
\bseries{DIMACS Series in Discrete Mathematics and Theoretical
Computer Science}
\bvolume{61}
\bpages{163--183}.
\bpublisher{Amer. Math. Soc.}, \baddress{Providence, RI}.
\bid{mr={1982426}}
\bptok{imsref}%
\end{bincollection}
%
\endbibitem

\bibitem{chak10}
%
\begin{bmisc}[author]
\bauthor{\bsnm{Chakerian},~\bfnm{J.}\binits{J.}} \AND
\bauthor{\bsnm{Holmes},~\bfnm{S.}\binits{S.}}
(\byear{2010}).
\bhowpublished{Computational tools for evaluating phylogenetic and hierachical
clustering trees.
Available at \href{http://arxiv.org/abs/1006.1015}{arXiv:1006.1015}.}
\bptok{imsref}%
\end{bmisc}
%
\endbibitem

\bibitem{deg05}
%
\begin{barticle}[pbm]
\bauthor{\bsnm{Degnan},~\bfnm{James~H.}\binits{J.~H.}} \AND
\bauthor{\bsnm{Salter},~\bfnm{Laura~A.}\binits{L.~A.}}
(\byear{2005}).
\btitle{Gene tree distributions under the coalescent process}.
\bjournal{Evolution}
\bvolume{59}
\bpages{24--37}.
\bid{issn={0014-3820}, pmid={15792224}}
\bptok{imsref}%
\end{barticle}
%
\endbibitem

\bibitem{don95}
%
\begin{barticle}[pbm]
\bauthor{\bsnm{Donnelly},~\bfnm{P.}\binits{P.}} \AND
\bauthor{\bsnm{Tavar{\'{e}}},~\bfnm{S.}\binits{S.}}
(\byear{1995}).
\btitle{Coalescents and genealogical structure under neutrality}.
\bjournal{Annu. Rev. Genet.}
\bvolume{29}
\bpages{401--421}.
\bid{doi={10.1146/annurev.ge.29.120195.002153}, issn={0066-4197},
pmid={8825481}}
\bptok{imsref}%
\end{barticle}
%
\endbibitem

\bibitem{doo99}
%
\begin{barticle}[author]
\bauthor{\bsnm{Doolittle},~\bfnm{W.~Ford}\binits{W.~F.}}
(\byear{1999}).
\btitle{Lateral genomics}.
\bjournal{Trends Genet.}
\bvolume{15}
\bpages{M5--M8}.
\bptok{imsref}%
\end{barticle}
%
\endbibitem

\bibitem{fel85}
%
\begin{barticle}[author]
\bauthor{\bsnm{Felsenstein},~\bfnm{J.}\binits{J.}}
(\byear{1985}).
\btitle{Confidence limits on phylogenies: An approach using the bootstrap}.
\bjournal{Evolution}
\bvolume{39}
\bpages{783--791}.
\bptok{imsref}%
\end{barticle}
%
\endbibitem

\bibitem{fel04}
%
\begin{bbook}[author]
\bauthor{\bsnm{Felsenstein},~\bfnm{J.}\binits{J.}}
(\byear{2004}).
\btitle{Inferring Phylogenies}.
\bpublisher{Sinauer}, \baddress{Sunderland, MA}.
\bptok{imsref}%
\end{bbook}
%
\endbibitem

\bibitem{flet04}
%
\begin{barticle}[author]
\bauthor{\bsnm{Fletcher},~\bfnm{P.~T}\binits{P.~T.}},
\bauthor{\bsnm{Lu},~\bfnm{C.}\binits{C.}},
\bauthor{\bsnm{Pizer},~\bfnm{S.~M.}\binits{S.~M.}} \AND
\bauthor{\bsnm{Joshi},~\bfnm{S.}\binits{S.}}
(\byear{2004}).
\btitle{Principal geodesic analysis for the study of nonlinear
statistics of
shape}.
\bjournal{IEEE Trans. Medical Imaging}
\bvolume{23}
\bpages{995--1005}.
\bptok{imsref}%
\end{barticle}
%
\endbibitem

\bibitem{good98}
%
\begin{barticle}[author]
\bauthor{\bsnm{Goodman},~\bfnm{M.}\binits{M.}},
\bauthor{\bsnm{Porter},~\bfnm{C.~A.}\binits{C.~A.}},
\bauthor{\bsnm{Czelusniak},~\bfnm{J.}\binits{J.}},
\bauthor{\bsnm{Page},~\bfnm{S.~L.}\binits{S.~L.}},
\bauthor{\bsnm{Schneider},~\bfnm{H.}\binits{H.}},
\bauthor{\bsnm{Shoshani},~\bfnm{J.}\binits{J.}},
\bauthor{\bsnm{Gunnell},~\bfnm{G.}\binits{G.}} \AND
\bauthor{\bsnm{Groves},~\bfnm{C.~P.}\binits{C.~P.}}
(\byear{1998}).
\btitle{Toward a phylogenetic classification of primates based on DNA evidence
complemented by fossil evidence}.
\bjournal{Mol. Phyl. Evol.}
\bvolume{9}
\bpages{585--598}.
\bptok{imsref}%
\end{barticle}
%
\endbibitem

\bibitem{gro87}
%
\begin{bincollection}[mr]
\bauthor{\bsnm{Gromov},~\bfnm{M.}\binits{M.}}
(\byear{1987}).
\btitle{Hyperbolic groups}.
In \bbooktitle{Essays in Group Theory}.
\bseries{Mathematical Sciences Research Institute Publications}
\bvolume{8}
\bpages{75--263}.
\bpublisher{Springer}, \baddress{New York}.
\bid{mr={0919829}}
\bptok{imsref}%
\end{bincollection}
%
\endbibitem

\bibitem{guin03}
%
\begin{barticle}[author]
\bauthor{\bsnm{Guindon},~\bfnm{S.}\binits{S.}} \AND
\bauthor{\bsnm{Gascuel},~\bfnm{O.}\binits{O.}}
(\byear{2003}).
\btitle{A simple, fast, and accurate algorithm to estimate large
phylogenies by
maximum likelihood}.
\bjournal{Syst. Biol.}
\bvolume{52}
\bpages{696--704}.
\bptok{imsref}%
\end{barticle}
%
\endbibitem

\bibitem{has89}
%
\begin{barticle}[mr]
\bauthor{\bsnm{Hastie},~\bfnm{Trevor}\binits{T.}} \AND
\bauthor{\bsnm{Stuetzle},~\bfnm{Werner}\binits{W.}}
(\byear{1989}).
\btitle{Principal curves}.
\bjournal{J. Amer. Statist. Assoc.}
\bvolume{84}
\bpages{502--516}.
\bid{issn={0162-1459}, mr={1010339}}
\bptok{imsref}%
\end{barticle}
%
\endbibitem

\bibitem{hill05}
%
\begin{barticle}[pbm]
\bauthor{\bsnm{Hillis},~\bfnm{David~M.}\binits{D.~M.}},
\bauthor{\bsnm{Heath},~\bfnm{Tracy~A.}\binits{T.~A.}} \AND
\bauthor{\bsnm{St.~John},~\bfnm{Katherine}\binits{K.}}
(\byear{2005}).
\btitle{Analysis and visualization of tree space}.
\bjournal{Syst. Biol.}
\bvolume{54}
\bpages{471--482}.
\bid{doi={10.1080/10635150590946961}, issn={1063-5157}, pii={L37747M5V37R177K},
pmid={16012112}}
\bptok{imsref}%
\end{barticle}
%
\endbibitem

\bibitem{hol05}
%
\begin{bincollection}[author]
\bauthor{\bsnm{Holmes},~\bfnm{S.}\binits{S.}}
(\byear{2005}).
\btitle{Statistical approach to tests involving phylogenies}.
In \bbooktitle{Mathematics of Evolution and Phylogeny}
(\beditor{\bfnm{O.}\binits{O.}~\bsnm{Gascuel}}, ed.)
\bpages{91--120}.
\bpublisher{Oxford Univ. Press}, \baddress{Oxford}.
\bptok{imsref}%
\end{bincollection}
%
\endbibitem

\bibitem{hue93}
%
\begin{barticle}[author]
\bauthor{\bsnm{Huelsenbeck},~\bfnm{J.}\binits{J.}} \AND
\bauthor{\bsnm{Hillis},~\bfnm{D.}\binits{D.}}
(\byear{1993}).
\btitle{Success of phylogenetic methods in the four-taxon case}.
\bjournal{Syst. Biol.}
\bvolume{42}
\bpages{247--264}.
\bptok{imsref}%
\end{barticle}
%
\endbibitem

\bibitem{kup08}
%
\begin{barticle}[mr]
\bauthor{\bsnm{Kupczok},~\bfnm{Anne}\binits{A.}},
\bauthor{\bsnm{Von~Haeseler},~\bfnm{Arndt}\binits{A.}} \AND
\bauthor{\bsnm{Klaere},~\bfnm{Steffen}\binits{S.}}
(\byear{2008}).
\btitle{An exact algorithm for the geodesic distance between phylogenetic
trees}.
\bjournal{J. Comput. Biol.}
\bvolume{15}
\bpages{577--591}.
\bid{issn={1066-5277}, mr={2425443}}
\bptok{imsref}%
\end{barticle}
%
\endbibitem

\bibitem{nye08}
%
\begin{barticle}[author]
\bauthor{\bsnm{Nye},~\bfnm{T.~M.~W.}\binits{T.~M.~W.}}
(\byear{2008}).
\btitle{Trees of trees: An approach to comparing multiple alternative
phylogenies}.
\bjournal{Syst. Biol.}
\bvolume{57}
\bpages{785--794}.
\bptok{imsref}%
\end{barticle}
%
\endbibitem

\bibitem{Nye11}
%
\begin{bmisc}[author]
\bauthor{\bsnm{Nye},~\bfnm{Tom~M.~W.}\binits{T.~M.~W.}}
(\byear{2011}).
\bhowpublished{Supplement to ``Principal components analysis in the
space of phylogenetic trees.''
\href{http://dx.doi.org/10.1214/11-AOS915SUPP}{DOI:10.1214/11-AOS915SUPP}.}
\bptok{imsref}%
\end{bmisc}
%
\endbibitem

\bibitem{owen09}
%
\begin{barticle}[author]
\bauthor{\bsnm{Owen},~\bfnm{Megan}\binits{M.}} \AND
\bauthor{\bsnm{Provan},~\bfnm{J.~Scott}\binits{J.~S.}}
(\byear{2010}).
\btitle{A fast algorithm for computing geodesic distances in tree space}.
\bjournal{IEEE/ACM Trans. Comp. Biol. and Bioinf.}
\bvolume{8}
\bpages{2--13}.
\bptok{imsref}%
\end{barticle}
%
\endbibitem

\bibitem{ram97}
%
\begin{barticle}[pbm]
\bauthor{\bsnm{Rambaut},~\bfnm{A.}\binits{A.}} \AND
\bauthor{\bsnm{Grassly},~\bfnm{N.~C.}\binits{N.~C.}}
(\byear{1997}).
\btitle{Seq-Gen: An application for the Monte Carlo simulation of DNA sequence
evolution along phylogenetic trees}.
\bjournal{Comput. Appl. Biosci.}
\bvolume{13}
\bpages{235--238}.
\bid{issn={0266-7061}, pmid={9183526}}
\bptok{imsref}%
\end{barticle}
%
\endbibitem

\bibitem{rob81}
%
\begin{barticle}[mr]
\bauthor{\bsnm{Robinson},~\bfnm{D.~F.}\binits{D.~F.}} \AND
\bauthor{\bsnm{Foulds},~\bfnm{L.~R.}\binits{L.~R.}}
(\byear{1981}).
\btitle{Comparison of phylogenetic trees}.
\bjournal{Math. Biosci.}
\bvolume{53}
\bpages{131--147}.
\bid{doi={10.1016/0025-5564(81)90043-2}, issn={0025-5564}, mr={0613619}}
\bptok{imsref}%
\end{barticle}
%
\endbibitem

\bibitem{sto02}
%
\begin{barticle}[pbm]
\bauthor{\bsnm{Stockham},~\bfnm{Cara}\binits{C.}},
\bauthor{\bsnm{Wang},~\bfnm{Li-San}\binits{L.-S.}} \AND
\bauthor{\bsnm{Warnow},~\bfnm{Tandy}\binits{T.}}
(\byear{2002}).
\btitle{Statistically based postprocessing of phylogenetic analysis by
clustering}.
\bjournal{Bioinformatics}
\bvolume{18}
\bpages{S285--S293}.
\bid{issn={1367-4803}, pmid={12169558}}
\bptok{imsref}%
\end{barticle}
%
\endbibitem

\bibitem{wang07}
%
\begin{barticle}[mr]
\bauthor{\bsnm{Wang},~\bfnm{Haonan}\binits{H.}} \AND
\bauthor{\bsnm{Marron},~\bfnm{J.~S.}\binits{J.~S.}}
(\byear{2007}).
\btitle{Object oriented data analysis: Sets of trees}.
\bjournal{Ann. Statist.}
\bvolume{35}
\bpages{1849--1873}.
\bid{doi={10.1214/009053607000000217}, issn={0090-5364}, mr={2363955}}
\bptok{imsref}%
\end{barticle}
%
\endbibitem

\end{thebibliography}
\end{document}